\documentclass[12pt]{amsart}

\usepackage{amsthm, amsmath, amssymb, mathrsfs}
\usepackage{enumerate}

\usepackage[colorlinks]{hyperref}

\makeatletter
\@namedef{subjclassname@2010}{%
	\textup{2010} Mathematics Subject Classification}
\makeatother

\allowdisplaybreaks

\theoremstyle{plain}
\newtheorem{thm}{Theorem}[section]
\newtheorem{cor}[thm]{Corollary}
\newtheorem{lem}[thm]{Lemma}
\newtheorem{prop}[thm]{Proposition}
\theoremstyle{definition}

\theoremstyle{remark}

\numberwithin{equation}{section}

\title[Amenability and Invariant subspaces of $PM_\Psi(G)$]{Amenability and Invariant subspaces of the algebra of pseudomeasures}
\author[A. Dabra]{Arvish Dabra$^\star$}
\address{Arvish Dabra,\newline\indent Department of Mathematics,\newline\indent Indian Institute of Technology Delhi,\newline\indent New Delhi - 110016, India.}
\email{arvishdabra3@gmail.com}
\author[N. S. Kumar]{N. Shravan Kumar}
\address{N. Shravan Kumar,\newline\indent Department of Mathematics,\newline\indent Indian Institute of Technology Delhi,\newline\indent New Delhi - 110016, India.}
\email{shravankumar.nageswaran@gmail.com}

\begin{document}
	
	\begin{abstract}
		Let $G$ be a locally compact group and $(\Phi,\Psi)$ a complimentary pair of Young functions. In this article, we consider the Banach algebra of $\Psi$-pseudomeasures $PM_\Psi(G)$ and the Orlicz Fig\`{a}-Talamanca Herz algebra $A_\Phi(G).$ We prove sufficient conditions for a group $G$ to be amenable in terms of the norm closed topologically invariant subspaces of $PM_\Psi(G).$ Further, for an amenable group $G$ with the Young function $\Phi$ satisfying the MA condition, we establish a one-to-one correspondence between certain topologically invariant subalgebras of $PM_\Psi(G)$ and the class of closed subgroups of $G.$ Moreover, we prove a similar result for the predual $A_\Phi(G)$ and derive a bijection between certain topologically invariant subalgebras of $A_\Phi(G)$ and the set of compact subgroups of $G.$
	\end{abstract}
	
	\keywords{Amenability; Banach algebras; Fig\`{a}-Talamanca Herz algebras; Orlicz spaces; Topologically invariant subspaces.\\
		$\star$ Corresponding author. \\
		{\it Email address:} arvishdabra3@gmail.com (A. Dabra).}
	
	\subjclass[2020]{Primary 43A15, 43A25, 46J10; Secondary 46E30, 46J99}
	
	\maketitle
	
	\section{Introduction} 

    For any locally compact group $G,$ the Group algebra $L^1(G),$ the Measure algebra $M(G),$ the Fourier algebra $A(G),$ the von Neumann algebra $VN(G)$ and certain convolution algebras are among the most widely studied and significant Banach algebras. Since the algebras $L^1(G)$ and $M(G)$ are \textit{non-commutative} Banach algebras in the case of \textit{non-abelian} groups, the analysis of their algebraic properties is limited. To overcome this limitation, Eymard \cite{Eym} pioneered the study of the Fourier algebra $A(G)$ and extensively studied the functorial properties of the \textit{commutative} Banach algebra $A(G)$ and its dual, $VN(G).$ If the underlying group $G$ is \textit{abelian}, then the Fourier algebra $A(G)$ is isometrically isomorphic to $L^1(\widehat{G}),$ where $\widehat{G}$ is the dual group of $G.$

    For $1 < p < \infty$ and any locally compact \textit{abelian} group $G,$ Fig\`{a}-Talamanca \cite{figa} introduced and studied the $L^p$ analogue of the Fourier algebra $A(G),$ denoted by $A_p(G).$ In 1971, Herz \cite{herz0} generalised the algebra $A_p(G)$ for \textit{any} locally compact group $G$ and termed as Fig\`{a}-Talamanca Herz algebra. For $p = 2,$ the algebra $A_p(G)$ coincides with the Fourier algebra $A(G).$ The algebra of $p$-pseudomeasures $PM_p(G)$ is the $w^\ast$-closure of $\lambda_p(M(G))$ in $\mathcal{B}(L^p(G)),$ where $\lambda_p$ denotes the representation of $M(G)$ on $L^p(G)$ given by left convolution. It is proved that the dual of $A_p(G)$ is isometrically isomorphic to $PM_{p'}(G),$ where $p$ and $p'$ are dual H\"{o}lder exponents.
    
    The algebra $A_p(G)$ has received significant attention from many researchers, such as Daws \cite{daws1}, Derighetti \cite{deriS, deri1, deriB}, Forrest \cite{FO1, FO2}, Granirer \cite{Gra0, Gra, GraS}, Hosseini and Amini \cite{Amini} and Lau and \"{U}lger \cite{lau} and others in the recent decades. Derighetti et al. \cite{DFM} studied the ideal structure of some Banach algebras related to the convolution operators on $L^p(G).$ The convolution algebras on $L^p$ spaces and the algebra of $p$-pseudomeasures $PM_p(G)$ have been recently studied by Daws and Spronk \cite{daws2} and Gardella and Thiel \cite{GH1, GH2}. In 2022, Gardella and Thiel \cite{GH2} proved that if $CV_p(G)$ is isometrically isomorphic to $CV_q(H),$ with $p,q \neq 2,$ then $G$ is (topologically) isomorphic to $H$ and $p$ and $q$ are either equal or conjugate. A similar result also holds for the algebra $PM_p(G).$
    
    It is well known that the Orlicz spaces are an important generalisation of the classical Lebesgue spaces and have garnered considerable interest from many authors in the last decade. \"{O}ztop et al. have intensively studied the weighted \cite{oztop1, oztop2} and twisted \cite{oztop3, oztop4} Orlicz algebras. In 2023, Tabatabaie and Latifpour \cite{tabata} studied some isomorphism preserving conditions for weighted Orlicz spaces. Further, the authors answered in affirmative to the problem of identifying groups that can be recovered from one of their convolution algebras. In particular, it is proved that for a class of Young functions $\Phi,$ if $CV_\Phi(G_1)$ and $CV_\Phi(G_2)$ are isometrically isomorphic, then $G_1$ is topologically isomorphic to $G_2,$ where $CV_\Phi(G)$ denotes the space of all convolution operators on the Orlicz space $L^\Phi(G),$ for any locally compact group $G.$

    Recently, Lal and Kumar \cite{RLSK1} and Aghababa and Akbarbaglu \cite{AA} independently introduced and studied the Orlicz $(L^\Phi)$ analogue of the Fig\`{a}-Talamanca Herz algebra $A_p(G),$ by replacing $L^p$ spaces with $L^\Phi,$ where $\Phi$ is a Young function satisfying the $\Delta_2$-condition. The $L^\Phi$ version of the algebra $A_p(G)$ is denoted by $A_\Phi(G)$ and termed as Orlicz Fig\`{a}-Talamanca Herz algebra. If $(\Phi,\Psi)$ corresponds to a complementary pair of Young functions, then the dual of $A_\Phi(G)$ is isometrically isomorphic to the space of $\Psi$-pseudomeasures $PM_\Psi(G)$ \cite[Theorem 3.5]{RLSK1}. The algebras $A_\Phi(G)$ and $PM_\Psi(G)$ have received a great deal of attention in recent years; see \cite{AA, AAR, AD, ArRaSh, RLSK1, RLSK3, RLSK2, RLSK4}. In 2024, the structure of maximal regular ideals and minimal ideals in the dual of certain topologically introverted subspaces of $PM_\Psi(G)$ is studied \cite{ArRaSh}.

    In this article, we are interested in investigating the connection between certain topologically invariant subspaces/subalgebras of the algebra of $\Psi$-pseudomeasures $PM_\Psi(G)$ and the underlying group $G.$ In Section \ref{sec2}, we recall the necessary notations and terminologies that are required in the sequel. In Section \ref{sec3}, inspired by Granirer \cite{GraS}, we prove sufficient conditions for a group $G$ to be \textit{amenable} in terms of the norm closed topologically invariant subspaces of $PM_\Psi(G).$ In Section \ref{sec4}, we derive a bijection between certain topologically invariant subalgebras of $PM_\Psi(G)$ and the set of closed subgroups of an \textit{amenable} group $G$ under the assumption that the Young function $\Phi$ satisfies the MA condition. This is Theorem \ref{thm4.3} of the article and serves as an Orlicz analogue of Lau and \"{U}lger's result (\cite[Theorem 8.15]{lau}) for the algebra $PM_p(G).$ A similar result also holds for the predual $A_\Phi(G).$  In Section \ref{sec5}, Theorem \ref{5.3} establishes a one-to-one correspondence between certain topologically invariant subalgebras of $A_\Phi(G)$ and the class of compact subgroups of $G.$

	\section{Preliminaries}\label{sec2}
    We begin this section by recalling the basic terminology related to the Orlicz spaces.

    A symmetric convex function $\Phi: \mathbb{R} \to [0,\infty]$ is said to be a \textit{Young function} if it satisfies $\Phi(0) = 0$ and $\lim\limits_{x \to \infty} \Phi(x) = +\infty.$ For each Young function $\Phi,$ a corresponding convex function $\Psi$ can be defined as follows:
    $$\Psi(y):= \sup{\{x\, |y|-\Phi(x):x\geq 0\}}, \hspace{1cm} y\in\mathbb{R}.$$
    The function $\Psi$ is also a Young function and is referred to as the complementary function to $\Phi.$ Further, the pair $(\Phi,\Psi)$ (and $(\Psi,\Phi)$) is called a complementary pair of Young functions (or a complementary Young pair). For $1 < p < \infty,$ the function $\Phi(x) = |x|^p/p$ serves as an example of a Young function, with $\Psi(y) = |y|^q/q$ as its complementary function, where $p$ and $q$ are dual H\"{o}lder exponents. 

    Let $G$ be a locally compact group with a fixed left Haar measure $dx.$ A Young function $\Phi$ is said to satisfy the \textit{$\Delta_{2}$-condition}, denoted by $\Phi\in\Delta_{2},$ if there exists a constant $k > 0$ and $x_{0} > 0$ such that the inequality $$\Phi(2x)\leq k \, \Phi(x)$$ holds for all $x\geq x_{0}$ whenever $G$ is compact and with $x_{0}=0$ if $G$ is non compact. The pair $(\Phi,\Psi)$ of complementary Young functions is said to satisfy the $\Delta_2$-condition if both $\Phi$ and $\Psi$ satisfy the $\Delta_2$-condition. It is easy to verify that the family of Young functions $\Phi_{\alpha,p}: y \mapsto \alpha \, |y|^p$ for $p \geq 1$ and $\alpha > 0,$ satisfies the $\Delta_2$-condition. Another example of a $\Delta_2$-function is the Young function $\Phi(x) = (1+|x|)\log(1+|x|)-|x|.$ However, its complementary function $\Psi(y) = e^{|y|}-|y|-1 \notin \Delta_2.$ 

    A Young function $\Phi$ is said to satisfy the \textit{Milnes-Akimoni$\check{c}$} condition (in short, MA condition) if for each $\epsilon > 0,$ there exists $\alpha_\epsilon > 1$ and an $x_1(\epsilon) \geq 0$ such that $$\Phi'((1+\epsilon)x) \geq \alpha_\epsilon \, \Phi'(x), \hspace{1cm} x \geq x_1(\epsilon).$$ The Young function $\Phi(x) = e^x - x - 1, \, 0 \leq x \in \mathbb{R},$ is an example of a Young function that satisfies the MA condition. Although its complementary function $\Psi(y) = (1+y)\log(1+y)-y, \, 0 \leq y \in \mathbb{R},$ does not satisfy the MA condition. It is easy to verify that the family of Young functions $\Phi_\alpha(x) = |x|^\alpha (1 + |\log|x||),$ $\alpha > 1,$ satisfies the MA condition. Moreover, the family $(\Phi_\alpha,\Psi_\alpha)$ of complementary pair of Young functions also satisfies the $\Delta_2$-condition.

    For a Young function $\Phi,$ the associated \textit{Orlicz space} is denoted by $L^{\Phi}(G)$ and is defined as follows:
    $${L}^{\Phi}(G) := \left\{ f: G \rightarrow  \mathbb{C}:f \, \mbox{is measurable and}\int_G\Phi(\alpha |f|)\ <\infty \text{ for some}~ \alpha>0  \right\}.$$ The Orlicz space $L^\Phi(G)$ is a Banach space when endowed with the Luxemburg (or Gauge) norm, given by
    $$N_{\Phi}(f):= \inf \left\{a>0:\int_G\Phi\left(\frac{|f|}{a}\right) \leq1 \right\}.$$ Let $\Psi$ be the complementary function to $\Phi.$ The Orlicz norm $\|\cdot\|_{\Phi}$ on $L^\Phi(G),$ is defined by $$\|f\|_{\Phi} := \sup \left\{\int_{G}|fg| :g\in L^\Psi(G) \, \, \text{and} \, \, \int_{G}\Psi(|g|) \leq1 \right\}.$$ In fact, these two norms are equivalent and for any $f \in L^\Phi(G),$ the following holds: $$N_\Phi(f) \leq \|f\|_\Phi \leq 2 \, N_\Phi(f).$$ 

    If a Young function $\Phi \in \Delta_{2},$ then the space $\mathcal{C}_c(G)$ of all continuous functions on $G$ with compact support is dense in $L^\Phi(G).$ Further, if the complementary function $\Psi$ is continuous and satisfies $\Psi(x) = 0$ if and only if $x = 0,$ then the dual space of $(L^\Phi(G),N_\Phi(\cdot))$ is isometrically isomorphic to $(L^\Psi(G),\|\cdot\|_\Psi)$ \cite[Corollary 9, Pg. 111]{RR}. In particular, if the complementary Young pair $(\Phi,\Psi)$ satisfies the $\Delta_2$-condition, then the Banach space $L^\Phi(G)$ is reflexive \cite[Theorem 10, Pg. 112]{RR}.
    
    We refer the readers to \cite{RR} for more details on the Orlicz spaces.
    \newline

    Let $\Phi$ be a Young function satisfying the $\Delta_2$-condition. For any function $g:G\rightarrow\mathbb{C}$, we define $\check{g}$ as $\check{g}(x):= g(x^{-1})$ for all $x \in G.$ The space $A_\Phi(G)$ is defined as the set of all continuous functions $u$ on $G$ that are of the form $$u=\sum_{n \in \mathbb{N}} f_n \ast \check{g_n},$$ where $f_n\in L^\Phi(G)$ and $g_n\in L^\Psi(G)$ such that $$\sum_{n \in \mathbb{N}} N_\Phi(f_n)\|g_n\|_\Psi<\infty.$$ Note that if $u \in A_\Phi(G),$ then $u \in \mathcal{C}_0(G).$ For any $u \in A_\Phi(G),$ the norm $\|\cdot\|_{A_\Phi}$ is given by $$\|u\|_{A_\Phi} := \inf\left\{\sum\limits_{n \in \mathbb{N}} N_\Phi(f_n) \|g_n\|_\Psi: u = \sum\limits_{n \in \mathbb{N}} f_n \ast \check{g_n}\right\}.$$ Equipped with this norm and pointwise addition and multiplication, the space $A_\Phi(G)$  becomes a commutative Banach algebra. The algebra $A_\Phi(G)$ is called the \textit{Orlicz Fig\`a-Talamanca Herz algebra}. Further, $A_\Phi(G)$ is a regular, tauberian, semi-simple Banach algebra with the Gelfand spectrum homeomorphic to $G.$

    Let $B_\Phi(G)$ denote the multiplier algebra of $A_\Phi(G)$ and defined by $$\{u \in \mathcal{C}(G): u \, v \in A_\Phi(G) \, \forall \, v \in A_\Phi(G)\}.$$ Each $u \in B_\Phi(G)$ defines a bounded linear map $M_u: A_\Phi(G) \to A_\Phi(G)$ given by $M_u(v) = u \, v,$ for $v \in A_\Phi(G).$ The space $B_\Phi(G)$ is a commutative Banach algebra endowed with the operator norm and pointwise addition and multiplication. Further, it is clear that $A_\Phi(G) \subseteq B_\Phi(G).$

    Let $M(G)$ be the set of all bounded complex Radon measures. For $\mu \in M(G)$ and $f \in L^\Psi(G),$ define the map $T_\mu: L^\Psi(G) \to L^\Psi(G)$ by left convolution, i.e., $T_\mu(f) := \mu \ast f.$ If $\mathcal{B}(L^\Psi(G))$ denotes the space of all bounded linear operators on $L^\Psi(G)$ with the operator norm, then it is easy to verify that $T_\mu \in \mathcal{B}(L^\Psi(G)).$ Let $PM_\Psi(G)$ be the closure of the set $\{T_\mu: \mu \in M(G)\}$ in $\mathcal{B}(L^\Psi(G))$ with respect to the ultra-weak topology ($w^*$-topology). It is proved in \cite{RLSK1} that for a locally compact group $G,$ the dual of $A_\Phi(G)$ is isometrically isomorphic to $PM_\Psi(G).$ Let $PF_\Psi(G)$ be the norm closure of the set $\{T_{\mu_f}: f \in L^1(G)\}$ inside $\mathcal{B}(L^\Psi(G)).$ The dual of $PF_\Psi(G)$ is $W_\Phi(G),$ where $W_\Phi(G)$ is a commutative Banach algebra containing $A_\Phi(G)$ \cite{RLSK4}. For $\mu = \delta_x \, \, (x \in G),$ we denote the operator $T_{\delta_x}$ by $\lambda_\Psi(x).$ Observe that $\lambda_\Psi$ represents the left regular representation of $G$ on $\mathcal{B}(L^\Psi(G)).$ From \cite{AD}, \cite{ArRaSh} and \cite{RLSK3}, let us recall the following subspaces of $PM_\Psi(G)$ that are used in this article.
    \begin{align*}
    C_{\delta,\Psi}(\widehat{G}) &:= \overline{span\{\lambda_\Psi(x):x \in G\}}, \\ M(\widehat{G}) &:= \overline{M(G)}, \\ UCB_\Psi(\widehat{G}) &:= \overline{A_\Phi(G)\cdot PM_\Psi(G)}, \\ AP_\Psi(\widehat{G}) &:= \{ T \in PM_\Psi(G): A_\Phi(G) \to PM_\Psi(G), u \mapsto u \cdot T \, \, \text{is compact}\}, \\ WAP_\Psi(\widehat{G}) &:= \{ T \in PM_\Psi(G): A_\Phi(G) \to PM_\Psi(G), u \mapsto u \cdot T \, \, \text{is weakly compact}\}.
    \end{align*}
    We remark that these subspaces are an Orlicz analogue of the subspaces of $PM_{p}(G)$ which are considered in \cite{Gra}.

    For more details on the Banach algebra $A_\Phi(G)$ and related spaces, we refer the readers to the series of papers \cite{AA, AAR, AD, ArRaSh, RLSK1, RLSK3, RLSK2, RLSK4}.
    \newline
    
    Let $\mathcal{A}$ be a Banach algebra and $X$ be a closed linear subspace of $\mathcal{A}^{\ast}.$ Then $X$ is said to be \textit{left topologically invariant} if $a \cdot \varphi \in X$ for all $a \in \mathcal{A}$ and $\varphi \in X,$ where $\langle b,a \cdot \varphi \rangle := \langle a \, b, \varphi \rangle$ for $b \in \mathcal{A}.$ Given such subspace $X$ of $\mathcal{A}^{\ast},$ one can define a continuous linear functional $m \odot \varphi$ on $\mathcal{A}$ by $$\langle a, m \odot \varphi \rangle := \langle a \cdot \varphi, m \rangle,$$ for $a \in \mathcal{A}, \varphi \in X$ and $m \in X^{\ast}.$ If $m \odot \varphi \in X$ for all $\varphi \in X$ and $m \in X^{\ast},$ then $X$ is termed as a \textit{left topologically introverted} subspace of $\mathcal{A}^{\ast}.$

    Similarly, we have the notion of \textit{right} topologically introverted subspace of $\mathcal{A}^{\ast}.$ However, in our case, since $\mathcal{A} = A_\Phi(G)$ is a commutative Banach algebra, both these notions coincide. It follows from \cite[Lemma 3.1]{ArRaSh} that the subspaces $C_{\delta,\Psi}(\widehat{G}), PF_\Psi(G), M(\widehat{G}), AP_\Psi(\widehat{G}), WAP_\Psi(\widehat{G}), UCB_\Psi(\widehat{G})$ and $PM_\Psi({G}),$ are norm closed topologically introverted subspaces of $PM_\Psi(G).$

    A locally compact group $G$ is said to be \textit{amenable} if it admits a left-invariant mean, i.e., there exists a positive linear functional $\Lambda$ on $L^\infty(G)$ of norm one such that $\langle L_x f, \Lambda \rangle = \langle f,\Lambda \rangle$ for all $f \in L^\infty(G)$ and $x \in G.$ Some examples of amenable groups are abelian groups, compact groups and solvable groups. However, the free group on two generators is not amenable. 

    Throughout this article, $G$ denotes a locally compact group with a fixed Haar measure $dx$ and the pair $(\Phi,\Psi)$ of complementary Young functions satisfy the $\Delta_2$-condition.
	
	\section{Amenability and Invariant subspaces of $PM_\Psi(G)$}\label{sec3}
	In this section, the main objective is to prove sufficient conditions for a group $G$ to be amenable in terms of the norm closed topologically invariant subspaces of $PM_\Psi(G).$ The approach and results are influenced by Granirer's work \cite{GraS}.
	
	For topologically invariant subspaces $X$ and $Y$ of $PM_\Psi(G),$ we denote the space of all continuous $A_\Phi$-module homomorphisms from $X$ to $Y$ by $\mathcal{M}(X,Y),$ i.e.,
	$$\mathcal{M}(X,Y) := \left\{T\in \mathcal{B}(X,Y): T(u\cdot\varphi) = u\cdot[T(\varphi)] \, \, \text{for all} \, \, u \in A_\Phi(G), \varphi \in X\right\}.$$ 
	
	Let us begin with the following lemma.
	
	\begin{lem}\label{s1p1}
		For a norm closed topologically invariant subspace $X$ of $PM_\Psi(G),$ consider the map $\eta: (\overline{A_\Phi\cdot X})^\ast \to \mathcal{M}(X,PM_\Psi)$ defined by $\eta(m) := m_L,$ where $m_L(\varphi) = m \odot \varphi$ for $\varphi \in X$ and $m \in (\overline{A_\Phi \cdot X})^\ast.$ Then $\eta$ is a norm decreasing injective map satisfying $$u \cdot [\eta(m)\varphi] = \eta(u \cdot m)\varphi$$ for all $u \in A_\Phi(G).$ Further, if $X$ is also topologically introverted, then $\eta$ is into $\mathcal{M}(X).$ 
	\end{lem}
	
	\begin{proof}
		For $\varphi \in X,$ observe that
		$$\|\eta(m)\varphi\| = \|m_L(\varphi)\| = \|m \odot \varphi\| \leq \|m\| \, \|\varphi\|$$
		which implies that $\|\eta(m)\| = \|m_L\| \leq \|m\|,$ i.e., $\eta$ is norm decreasing.
		
		To prove that $\eta$ is injective, for $m \in (\overline{A_\Phi\cdot X})^\ast,$ let
		\begin{align*}
			\eta(m) = 0 &\iff m_L = 0\\
			&\iff m \odot \varphi = 0 \, \, \forall \, \, \varphi \in X\\
			&\iff \langle u, m \odot \varphi \rangle = 0 \, \, \forall \, \, u \in A_\Phi(G), \varphi \in X\\
			&\iff \langle u \cdot \varphi, m \rangle = 0 \, \, \forall \, \, u \in A_\Phi(G), \varphi \in X\\
			&\iff m = 0.
		\end{align*}
		
		Also, for $u \in A_\Phi(G), \varphi \in X$ and $m \in (\overline{A_\Phi \cdot X})^\ast,$
		$$u \cdot [\eta(m)\varphi] = u \cdot (m \odot \varphi) = (u \cdot m) \odot \varphi = \eta(u \cdot m)\varphi.$$
		
		If $X$ is topologically introverted, then $m \odot \varphi \in X$ for all $m \in (\overline{A_\Phi \cdot X})^\ast$ and $\varphi \in X.$ Hence, the result follows.
	\end{proof}
	
	The subsequent theorem provides that the map $\eta$ is an isomorphism under the assumption that the underlying group $G$ is amenable.
	
	\begin{thm}\label{S1p2}
		Let $G$ be an amenable group and $X$ be a norm closed topologically invariant subspace of $PM_\Psi(G).$ Then the spaces $(\overline{A_\Phi \cdot X})^\ast$ and $\mathcal{M}(X,PM_\Psi)$ are isomorphic with equivalent norms. If $X$ is also topologically introverted, then $\eta[(\overline{A_\Phi \cdot X})^\ast] = \mathcal{M}(X).$ Further, if $X \subseteq UCB_\Psi(\widehat{G}),$ then $\eta(X^\ast) = \mathcal{M}(X).$
	\end{thm}
	
	\begin{proof}
		Since $G$ is amenable, it follows from \cite[Theorem 3.1]{RLSK2} that $A_\Phi(G)$ has a bounded approximate identity (BAI), say $\{v_\alpha\}_{\alpha \in \Lambda},$ such that $\|v_\alpha\| \leq 2 \, \, \forall \, \, \alpha \in \Lambda.$
		
		For $\varphi \in X$ and $m \in (\overline{A_\Phi \cdot X})^\ast,$ observe that
		$$| \langle v_\alpha \cdot \varphi, m \rangle| = | \langle v_\alpha, m \odot \varphi \rangle| = |\langle v_\alpha, m_L(\varphi) \rangle| \leq \|m_L(\varphi)\| \, \|v_\alpha\| \leq 2 \, \|m_L(\varphi)\|.$$
		As $\{v_\alpha\}_{\alpha \in \Lambda}$ is a BAI, we have,
		\begin{equation}\label{s1e1}
			|\langle \varphi, m \rangle| \leq 2\, \|m_L(\varphi)\| \, \, \text{provided} \, \, \varphi \in (\overline{A_\Phi \cdot X}).
		\end{equation}
		
		Let $\epsilon > 0.$ Then there exists $\varphi_0 \in (A_\Phi \cdot X)$ with $\|\varphi_0\| \leq 1$ such that $$\|m\| - \epsilon \leq |\langle \varphi_0, m \rangle|.$$
		As $\varphi_0 \in (A_\Phi \cdot X),$ it follows from (\ref{s1e1}) that
		$$\|m\| - \epsilon \leq 2 \, \|m_L(\varphi_0)\| \leq 2 \, \|m_L\|.$$
		Thus, by Lemma \ref{s1p1}, we have,
		$$\|m_L\| \leq \|m\| \leq 2 \, \|m_L\|$$
		which proves that the norms are equivalent.
		
		It is easy to verify that $\eta$ is a homomorphism between $(\overline{A_\Phi \cdot X})^\ast$ and $\mathcal{M}(X,PM_\Psi).$ The crucial step is to prove the surjectivity of the map $\eta.$ Let $T \in \mathcal{M}(X,PM_\Psi)$ and $T^t$ be its transpose conjugate. Using the canonical identification of $A_\Phi(G)$ in its double dual $PM_\Psi(G)^*,$ we have,
		$$\|T^t(v_\alpha)\| \leq \|T^t\| \, \|v_\alpha\| \leq 2 \, \|T^t\|.$$
		Now, by Banach-Alaoglu Theorem, the net $\{T^t(v_\alpha)\}$ has a subnet say $\{T^t(v_{\alpha_\beta})\}$ which converges in $w^\ast$-topology. Let $m$ be the $w^\ast$-limit of this subnet in $(\overline{A_\Phi \cdot X})^\ast.$ Observe that for any $u \in A_\Phi(G)$ and $\varphi \in X,$
		\begin{align*}
			\langle u, T(\varphi) \rangle = \lim_\beta \, \langle u \, v_{\alpha_\beta}, T(\varphi) \rangle &= \lim_\beta \, \langle v_{\alpha_\beta}, u \cdot T(\varphi) \rangle \\ &= \lim_\beta \, \langle v_{\alpha_\beta}, T(u \cdot \varphi) \rangle \\ &= \lim_\beta \, \langle u \cdot \varphi, T^t(v_{\alpha_\beta}) \rangle \\ &= \langle u \cdot \varphi, m \rangle \\ &= \langle u, m \odot \varphi \rangle = \langle u, m_L(\varphi) \rangle.
		\end{align*}
		Thus, it follows that $T = m_L.$ This proves that the map $\eta$ is an isomorphism between $(\overline{A_\Phi \cdot X})^\ast$ and $\mathcal{M}(X,PM_\Psi)$ with equivalent norms.
		
		If $X$ is topologically introverted, then $m \odot \varphi \in X$ for all $m \in (\overline{A_\Phi \cdot X})^\ast$ and $\varphi \in X.$ Thus, $\eta[(\overline{A_\Phi \cdot X})^\ast] = \mathcal{M}(X).$
		
		Further, if $X \subseteq UCB_\Psi(\widehat{G}),$ then $(A_\Phi \cdot X) = X$ and hence, $\eta(X^\ast) = \mathcal{M}(X).$
	\end{proof}	
	
	The following corollary is a straightforward consequence of Theorem \ref{S1p2}.
	
	\begin{cor}\label{s1c1}
		For an amenable group $G,$ the spaces $UCB_\Psi(\widehat{G})^\ast$ and $\mathcal{M}(PM_\Psi)$ are isomorphic with equivalent norms. 
	\end{cor}
	
	\begin{cor}
		Let $G$ be an amenable group. Then the space $W_\Phi(G)$ is continuously embedded in $UCB_\Psi(\widehat{G})^\ast.$
	\end{cor}
	
	\begin{proof}
		For $w \in W_\Phi(G),$ define the map $\nu: W_\Phi(G) \to \mathcal{M}(PM_\Psi)$ by $$[\nu(w)](T) := w \cdot T \hspace{1cm} (T \in PM_\Psi(G)).$$ The map $\nu$ is well-defined as $PM_\Psi(G)$ is a $B_\Phi(G)$-module and $W_\Phi(G) \subseteq B_\Phi(G)$ \cite[Corollary 3.2 (ii)]{RLSK4}. This provides us with a continuous embedding of $W_\Phi(G)$ in $\mathcal{M}(PM_\Psi)$ and the final result is just a consequence of Corollary \ref{s1c1}.
	\end{proof}

	Here is the promised result of this section, which provides sufficient conditions for a group to be amenable.
	
	\begin{thm}
		Let $(\Phi,\Psi)$ satisfies the MA condition and $X$ be a norm closed topologically invariant subspace of $PM_\Psi(G)$ that contains $M_1^+(G),$ where $M_1^+(G) := \{0 \leq \mu \in L^1(G)^\ast:\mu(G) =1\},$ the set of positive probability Radon measures. If $\|T_\mu\|_{PM_\Psi} \leq 1$ for all $\mu \in M_1^+(G)$ and  $\eta[(\overline{A_\Phi \cdot X})^\ast]$ contains the identity injection $I \in \mathcal{M}(X,PM_\Psi)$ (a fortiori if $\eta$ is surjective), then $G$ is amenable.   
	\end{thm}
	
	Recall that $M(\widehat{G}), UCB_\Psi(\widehat{G}),        
	WAP_\Psi(\widehat{G})$ and $PM_\Psi(G)$ are some examples of such      $X.$ Further, we remark that if we equip the Orlicz space      $L^\Phi(G)$ with the Orlicz norm $\|\cdot\|_\Phi,$ then the dual space is given by $L^\Psi(G)$ with the Gauge norm $N_\Psi(\cdot)$ and as a consequence, $\|T_\mu\|_{PM_\Psi} \leq 1$ for all $\mu \in M(G)$ and thus, satisfies the hypothesis of the theorem.
	
	\begin{proof}
		Let $m \in (\overline{A_\Phi \cdot X})^\ast$ be such that $\eta(m) = I.$ By definition of $\eta,$ for $u \in A_\Phi(G)$ and $T_\mu \in X,$ where $\mu$ is any bounded Radon measure in $X,$ we have,
		$$\langle u \cdot T_\mu, m \rangle = \langle u, m \odot T_\mu \rangle = \langle u, \eta(m)T_\mu \rangle = \langle u, I(T_\mu) \rangle = \langle u, T_\mu \rangle = \int_G u \, d\mu.$$
		Thus, 
		\begin{equation}\label{s1e2}
			\left|\int_G u \, d\mu\right| = |\langle u \cdot T_\mu, m \rangle| \leq \|m\| \, \|u \cdot T_\mu\|_{PM_\Psi}.
		\end{equation}
		
		Given $\epsilon > 0.$ Choose $G_0$ in $G$ such that $G_0$ is $\sigma$-compact with $\mu(G \setminus G_0) < \epsilon.$ Further, for each $n \in \mathbb{N},$ one can choose compact sets $K_n$ such that $K_n \subseteq K_{n+1}$ and $\bigcup_{n \in \mathbb{N}} K_n = G_0.$ Now, choose $u_n \in A_\Phi(G) \cap C_c(G)$ such that $u_n(K_n) = 1$ and $0 \leq u_n \leq 1.$ Then, 
		\begin{equation}\label{s1e3}
			\|u_n \cdot T_\mu - T_\mu\|_{PM_\Psi} \leq \|u_n \cdot \mu - \mu\|_{M(G)} \to  0.
		\end{equation}
		By replacing $u$ with $u_n$ in (\ref{s1e2}) and using (\ref{s1e3}), it follows that
		$$|\mu(G)| \leq \|m\| \, \|T_\mu\|_{PM_\Psi}.$$
		As $X$ is closed under convolution powers, for all $n \in \mathbb{N},$
		$$|\mu(G)|^n = |\mu^{\ast^n}(G)| \leq \|m\| \, \|T_{\mu^{\ast^n}}\|_{PM_\Psi} = \|m\| \, \|T_\mu^n\|_{PM_\Psi}.$$
		This implies that for each $\mu \in M_1^+(G) \subseteq X,$
		$$1 = |\mu(G)| \leq \|m\|^{1/n} \, \|T_\mu^n\|_{PM_\Psi}^{1/n} \to \rho(T_\mu) \leq \|T_\mu\|_{PM_\Psi} \leq 1,$$
		i.e., the spectral radius $\rho$ of the operator $T_\mu$ is $1$ for each $\mu \in M_1^+(G).$ Hence, it follows from \cite[Theorem 8]{raoIII} that $G$ is amenable.
	\end{proof}
	
	\section{Invariant subalgebras of $PM_\Psi(G)$ and closed subgroups of $G$}\label{sec4}
	In this section, we establish a bijection between the set of closed subgroups of an amenable group $G$ and certain topologically invariant subalgebras of $PM_\Psi(G),$ assuming that the Young function $\Phi$ satisfies the MA condition. The results of this section are motivated by \cite{lau}. To start, we shall define the necessary notations.
	
	For a closed subset $K$ of $G,$ we denote the $w^\ast$-closure of the linear span of $\{\lambda_\Psi(x):x\in K\}$ in $PM_\Psi(G)$ by $C_{\Psi,K}(\widehat{G}).$ Without any ambiguity, we shall use the notation $C_K$ in place of $C_{\Psi,K}(\widehat{G}).$ It is easy to verify that $C_K$ is a $w^\ast$-closed topologically invariant subspace of $PM_\Psi(G).$ 
	
	For each $w^\ast$-closed topologically invariant subspace $C$ of $PM_\Psi(G),$ let $\Gamma(C)$ denote the closed subset $\{x \in G:\lambda_\Psi(x) \in C\}$ of $G.$ Note that for any closed subgroup $K$ of $G,$ it holds that $K \subseteq \Gamma(C_K).$ For each $T \in PM_\Psi(G),$ let $\langle T\rangle$ denote the $w^\ast$-closure of the linear span of $\{u\cdot T:u \in A_\Phi(G)\}.$ Observe that $\langle T\rangle$ is the smallest $w^\ast$-closed topologically invariant subspace of $PM_\Psi(G)$ that contains the set $\{u\cdot T:u \in A_\Phi(G)\}.$ As mentioned in the previous section, for an amenable group $G,$ the algebra $A_\Phi(G)$ has a bounded approximate identity, which implies that $T \in \langle T\rangle.$
	
	Recall from \cite[Pg. 101]{herz} that for an algebra of functions $\mathcal{A},$ the support of a linear functional $T \in \mathcal{A}^\ast$ as a subset of $G$ is characterised as follows: $x \notin supp(T)$ if and only if there exists a neighbourhood $W$ of $x$ such that $\langle w,T \rangle=0$ for all $w \in \mathcal{A}$ with $supp(w) \subseteq W.$
	
	We have the following straightforward yet significant lemma.
	
	\begin{lem}\label{s2l1}
		For any $T \in PM_\Psi(G),$ the support of $T$ is $\Gamma(\langle T\rangle).$ 
	\end{lem}
	
	\begin{proof}
		Let $x \in supp(T).$ Then $\lambda_\Psi(x)$ belongs to the $w^\ast$-closure of the $span\{u \cdot T: u \in A_\Phi(G)\}$ in $PM_\Psi(G).$ Thus, $\lambda_\Psi(x) \in \langle T\rangle$ and hence, $x \in \Gamma(\langle T\rangle)$ proving that $supp(T) \subseteq \Gamma(\langle T\rangle).$ To prove the reverse inclusion, by using the above characterisation of $supp(T),$ one can observe that $x \notin supp(T)$ implies that $\lambda_\Psi(x) \notin \langle T\rangle$ and hence, $x \notin \Gamma(\langle T\rangle).$
	\end{proof}
	
	\begin{prop}\label{s2l2}
		Let $H$ be a closed subgroup of an amenable group $G$ and $\Phi$ satisfies the MA condition. Then for any $T \in PM_\Psi(G),$ $supp(T) \subseteq H$ if and only if $T \in C_H.$
	\end{prop}
	
	\begin{proof}
		To prove that $T \in C_H$ implies $supp(T) \subseteq H,$ it is equivalent to prove that if $supp(T) \not\subseteq H$ then $T \notin C_H.$ Assume that $supp(T) \not\subseteq H$ and let $x \in supp(T) \setminus H.$ Since $x \in supp(T)$ and $G\setminus H$ is open, there exists a neighbourhood $U$ of $x,$ such that $\langle u,T\rangle \neq 0$ with $supp(u) \subseteq U \subseteq G \setminus H.$ This forces that $T \notin C_H$ as the $supp(u)$ is disjoint from $H.$
		
		Conversely, let  $supp(T) \subseteq H.$ By \cite[Corollary 3.4]{RLSK2}, since $H$ is a closed subgroup of an amenable group $G,$ it is a set of spectral synthesis for $A_\Phi(G).$ Now, consider the following closed ideal of $A_\Phi(G),$ given by $$I_\Phi(H) := \{u\in A_\Phi(G): u(H) =\{0\}\}.$$ By Lemma \ref{s2l1}, since $\Gamma(\langle T \rangle) = supp(T) \subseteq H,$ it follows from \cite[Theorem 40.10, Pg. 526]{ross} that
		$$T \in I_\Phi(H)^\perp = \left(^\perp\{\lambda_\Psi(x):x\in H\}\right)^\perp = \overline{span\{\lambda_\Psi(x):x\in H\}}^{w^\ast}.$$
		Thus, $T$ belongs to the $w^\ast$-closure of the linear span of $\{\lambda_\Psi(x):x\in H\},$ i.e., $T \in C_H.$ Hence, the result follows.
	\end{proof}
	
	Before proceeding to the main result of this section, we need to introduce a few additional notations. Let $\mathcal{H}$ be the collection of all closed subgroups of $G$ and let $\mathcal{C}$ represents the collection of all $w^\ast$-closed topologically invariant subalgebras of $PM_\Psi(G)$ such that, for $C \in \mathcal{C},$ $\Gamma(C)$ is a \textit{subgroup} of $G.$
	
	Here is the main theorem of this section, which serves as an Orlicz analogue of \cite[Theorem 8.15]{lau}.
	
	\begin{thm}\label{thm4.3}
		Let $G$ be an amenable group and $\Phi$ satisfies the MA condition. Then the mapping $\gamma: \mathcal{H} \to \mathcal{C}$ defined by $\gamma(H) = C_H$ is a bijection.
	\end{thm}
	
	\begin{proof}
		Let $H \in \mathcal{H}.$ By definition, $C_H$ is a $w^\ast$-closed topologically invariant subalgebra of $PM_\Psi(G).$ By Proposition \ref{s2l2}, it follows that $\Gamma(C_H) = H.$ Thus, the map $\gamma$ is well-defined.
		
		To prove that $\gamma$ is injective, let $H_1, H_2 \in \mathcal{H}$ be such that $H_1 \neq H_2,$ i.e., there exists $h \in H_1 \setminus H_2.$ Since a locally compact Hausdorff space is regular, there exists disjoint open sets $U$ and $V$ such that $h \in U$ and $H_2 \subseteq V.$ Choose $u \in A_\Phi(G) \cap C_c(G)$ such that $\langle u, \lambda_\Psi(h)\rangle = u(h) = 1$ and $u(H_2) = \{0\}.$ Then, $u \in I_\Phi(H_2)$ but $u \notin I_\Phi(H_1)$ which implies that $\lambda_\Psi(h) \in C_{H_1} \setminus C_{H_2}$ and hence, $C_{H_1} \neq C_{H_2}.$
		
		For surjectivity, let $C \in \mathcal{C}.$ Then, $\Gamma(C)$ is a closed subgroup of $G,$ say $H$ and clearly, $C_H \subseteq C.$ If $T \in C,$ by Lemma \ref{s2l1}, $$supp(T) = \Gamma(\langle T \rangle) \subseteq \Gamma(C) = H.$$
		Thus, it follows from Proposition \ref{s2l2} that $T \in C_H$ and hence, $C = C_H.$ 
	\end{proof}
	
	\begin{cor}
	Let $\{H_\alpha\}_{\alpha \in \Lambda}$ be a family of closed subgroups of an amenable group $G$ and $\Phi$ satisfies the MA condition. Then $$\gamma\left(\bigcap_{\alpha \in \Lambda} H_\alpha \right) = \bigcap_{\alpha \in \Lambda} \gamma(H_\alpha).$$
	\end{cor}
	
	\begin{proof}
		It is clear from the definition of $\gamma$ that if $H_1 \subseteq H_2$ then $\gamma(H_1) \subseteq \gamma(H_2).$ Thus,
		$$\gamma\left(\bigcap_{\alpha \in \Lambda} H_\alpha \right) \subseteq \bigcap_{\alpha \in \Lambda} \gamma(H_\alpha).$$
		
		The crucial step is to prove the reverse inclusion. If $T \in \bigcap\limits_{\alpha \in \Lambda} \gamma(H_\alpha),$ then $T \in \gamma(H_\alpha) = C_{H_\alpha}$ for each $\alpha \in \Lambda.$ Now, by Proposition \ref{s2l2}, it follows that $supp(T) \subseteq H_\alpha$ for each $\alpha \in \Lambda$ and thus, $supp(T) \subseteq \bigcap\limits_{\alpha \in \Lambda} H_\alpha.$ Again, by Proposition \ref{s2l2}, $$T \in C_{\bigcap\limits_{\alpha \in \Lambda} H_\alpha} = \gamma(\bigcap\limits_{\alpha \in \Lambda} H_\alpha).$$
		Hence, the result follows.
	\end{proof}
	
	\section{Invariant subalgebras of $A_\Phi(G)$ and compact subgroups of $G$}\label{sec5}
    
    In this section, we derive a one-to-one correspondence between certain topologically invariant subalgebras of $A_\Phi(G)$ and the set of compact subgroups of $G.$ The ideas are inspired by \cite{lau}. Let us begin by recalling some standard notations and results on the algebra $A_\Phi(G).$

    For $x \in G$ and $f \in L^\Phi(G),$ let $L_x$ and $R_x$ denote the left and right translation by $x$ on $L^\Phi(G),$ respectively. The definitions are as follows:
    $$(L_xf)(y):= f(x^{-1}y) \hspace{0.5cm} \text{and} \hspace{0.5cm} (R_xf)(y) := f(yx), \hspace{1cm} (y \in G).$$

    Recall from \cite[Definition 5.1]{AA} and \cite[Theorem 3.5]{AAR} that for any $x \in G$ and $f \in L^\Phi(G),$
    $$N_\Phi(L_xf) = N_\Phi(f) \hspace{0.5cm} \text{and} \hspace{0.5cm} \|L_x f\|_\Phi = \|f\|_\Phi.$$
    It is easy to verify that the algebra $A_\Phi(G)$ is left and right translation invariant and for any $x \in G$ and $u \in A_\Phi(G),$
    $$\|L_x u\| = \|u\| \hspace{0.5cm} \text{and} \hspace{0.5cm} \|R_x u\| = \|u\|,$$
    holds.
    
    Furthermore, as mentioned in \cite[Theorem 3.5]{AAR}, it follows from \cite[Proposition 5.3]{AA} that for any fixed $u \in A_\Phi(G),$ the mappings 
    $$x \mapsto L_x u \hspace{0.5cm} \text{and} \hspace{0.5cm} x \mapsto R_x u$$
    from $G$ to $A_\Phi(G)$ are continuous.

    It is well known that 
    $$A_\Phi(G) \subseteq \mathcal{C}_0(G) \subseteq L^\infty(G) = {L^1(G)}^\ast.$$
    So, one can equip $A_\Phi(G)$ with the $w^\ast$-topology induced on $A_\Phi(G)$ by $w^\ast$-topology on $L^\infty(G)$ and we denote this topology by $\sigma(A_\Phi(G),L^1(G)).$

    Now, for any subgroup $H$ of $G,$ consider the following subspace of $A_\Phi(G),$
    $$\mathcal{L}_\Phi(H) := \{ u \in A_\Phi(G): L_h u = u \, \, \forall \, \, h \in H \}.$$

    The subsequent lemma outlines the significant properties of the space $\mathcal{L}_\Phi(H)$ for any closed subgroup $H$ of $G.$

    \begin{lem}\label{5.1}
        For any closed subgroup $H$ of $G,$ the space $\mathcal{L}_\Phi(H)$ is a $\sigma(A_\Phi(G),L^1(G))$-closed right translation invariant subalgebra of $A_\Phi(G)$ that is closed under conjugation.
    \end{lem}

    \begin{proof}
        It is easy to verify that for any closed subgroup $H$ of $G,$ the space $\mathcal{L}_\Phi(H)$ is a subalgebra of $A_\Phi(G)$ that is closed under conjugation. Further, for any $u \in \mathcal{L}_\Phi(H)$ and $g \in G,$ $$L_h(R_g u) = R_g(L_h u) = R_g u,$$
        for all $h \in H.$ This implies that $R_g u \in \mathcal{L}_\Phi(H)$ for all $g \in G,$ i.e., $\mathcal{L}_\Phi(H)$ is right translation invariant. 

        Further, observe that 
        $$\mathcal{L}_\Phi(H) = \bigcap\limits_{h \in H} \{u \in A_\Phi(G) : L_h u - u = 0\}.$$
        Now, for each $h \in H,$ consider the map $f \mapsto L_h f - f$ from $L^\infty(G)$ to itself. By restricting this map to $A_\Phi(G)$ and using the $w^\ast$-$w^\ast$-continuity, it follows that 
        $$\{u \in A_\Phi(G) : L_h u - u = 0\}$$
        is $\sigma(A_\Phi(G),L^1(G))$-closed for each $h \in H.$ Since arbitrary intersection of closed sets is closed, the subalgebra $\mathcal{L}_\Phi(H)$ is also $\sigma(A_\Phi(G),L^1(G))$-closed.
    \end{proof}
        
        The next proposition serves as a pivotal step in proving the main theorem of this section and gives that $\mathcal{L}_\Phi(H) \neq \{0\}$ for any compact subgroup $H$ of $G.$

    \begin{prop}\label{5.2}
        If $H$ is a compact subgroup of $G,$ then the subalgebra $\mathcal{L}_\Phi(H)$ is non-zero. Further, for $g \notin H,$ there exists $\Tilde{u}_g \in \mathcal{L}_\Phi(H)$ such that $L_g \Tilde{u}_g \neq \Tilde{u}_g.$
    \end{prop}

    \begin{proof}
        As $H$ is compact, choose $u \in A_\Phi(G)$ such that $u(H) = \{1\}.$ Since the map $g \mapsto L_g u$ is continuous from $G$ to $A_\Phi(G),$ the set $\{L_h u : h \in H\}$ is compact in $A_\Phi(G).$ Let $\Bar{E}_\Phi(u,H)$ denote the closed convex hull of this set in $A_\Phi(G)$ and it follows from \cite[Theorem 5.35]{guide} that $\Bar{E}_\Phi(u,H)$ is also compact. Let $\Sigma$ be the group $\{L_h : h \in H\}$ of continuous affine transformations from $\Bar{E}_\Phi(u,H)$ to itself. Since the map $h \mapsto L_h$ from $H$ to $\Sigma$ is continuous, it follows that $\Sigma$ is also a compact semigroup. Thus, by \cite[Corollary 1]{rosen}, there exists an invariant mean on $\Sigma$ and therefore implies that $\Sigma$ is amenable. Now, by \cite[Theorem 1]{day}, there exists a common fixed point of $\Sigma$ in $\Bar{E}_\Phi(u,H),$ say $\Tilde{u},$ i.e.,
        $$L_h \Tilde{u} = \Tilde{u} \, \, \forall \, \, h \in H,$$
        which implies that $\Tilde{u} \in \mathcal{L}_\Phi(H).$ Since $\Tilde{u} \in \Bar{E}_\Phi(u,H)$ and $L_hu(e) = u(h^{-1}) = 1$ for any $h \in H,$ it follows that $\Tilde{u}(e) = 1,$ thus proving that $\mathcal{L}_\Phi(H) \neq \{0\}.$ Further, observe that 
        $$\Tilde{u}(h) = L_{h^{-1}}\Tilde{u}(e) = \Tilde{u}(e) = 1 \, \, \forall \, \,  h \in H.$$
        Now, for $g \notin H,$ $Hg^{-1}$ and $H$ are disjoint compact subsets of $G.$ Since $A_\Phi(G)$ is a tauberian algebra, there exists $u_g \in A_\Phi(G)$ such that $u_g(H) = \{1\}$ and $u_g(Hg^{-1}) = \{0\}.$ By repeating the same argument for $u_g$ in place of $u,$ one can find $\Tilde{u}_g$ in the closed convex hull of $\{L_h u_g : h \in H\}$ such that $\Tilde{u}_g \in \mathcal{L}_\Phi(H)$ and $\Tilde{u}_g(e) = 1.$ Since $L_h u_g(g^{-1}) = u_g(h^{-1}g^{-1}) = 0,$ we have, $\Tilde{u}_g(g^{-1}) = 0.$ Thus, 
        $$L_g \Tilde{u}_g(e) = \Tilde{u}_g(g^{-1}) = 0 \neq 1 = \Tilde{u}_g(e),$$
        and hence, the result follows.
    \end{proof}

    Here is the promised result of this section, which provides a bijection between certain topologically invariant subalgebras of $A_\Phi(G)$ and the class of compact subgroups of $G.$ This theorem is an Orlicz analogue of \cite[Theorem 8.11]{lau} and a weak generalisation of \cite[Theorem 9]{tak}.

    \begin{thm}\label{5.3}
        The map $\mathcal{F}$ given by $$\mathcal{F}(H) = \mathcal{L}_\Phi(H)$$
        is a bijection between the set of compact subgroups of $G$ and the set of all non-zero $\sigma(A_\Phi(G),L^1(G))$-closed right translation invariant subalgebras of $A_\Phi(G)$ that are closed under conjugation.
    \end{thm}

    \begin{proof}
        It follows from Lemma \ref{5.1} and Proposition \ref{5.2} that the map $\mathcal{F}$ is well-defined. To prove that $\mathcal{F}$ is one-one, let $H_1$ and $H_2$ be two distinct compact subgroups of $G.$ Then there exists $g \in H_2 \setminus H_1$ and by Proposition \ref{5.2}, there exists $\Tilde{u}_{g} \in \mathcal{L}_\Phi(H_1)$ such that $L_{g} \Tilde{u}_{g} \neq \Tilde{u}_{g}.$ This implies that $\Tilde{u}_{g} \notin \mathcal{L}_\Phi(H_2)$ and hence, $\mathcal{L}_\Phi(H_1) \neq \mathcal{L}_\Phi(H_2).$

        The crucial and challenging step is to prove that the map $\mathcal{F}$ is surjective. Let $\mathbf{L}$ be a non-zero $\sigma(A_\Phi(G),L^1(G))$-closed right translation invariant subalgebra of $A_\Phi(G)$ which is closed under conjugation. Consider the following subset $G_\mathbf{L}$ of $G$ given by 
        $$G_\mathbf{L} := \{g \in G: L_g u = u \, \, \forall \, \, u \in \mathbf{L}\}.$$
        It is easy to verify that $G_\mathbf{L}$ is indeed a closed subgroup of $G.$ We claim that $G_\mathbf{L}$ is also compact. Let $\theta(A_\Phi(G))$ be the Gelfand spectrum of $A_\Phi(G)$ and $\kappa$ be the map from $G$ to $\theta(A_\Phi(G)),$ given by $\kappa(g) := \lambda_\Psi(g),$ where $\langle u,\lambda_\Psi(g) \rangle := u(g),$ for $u \in A_\Phi(G).$ By \cite[Corollary 3.8]{RLSK1}, it follows that $\kappa$ is a homeomorphism. Since $\theta(A_\Phi(G)) \cup \{0\}$ is $w^\ast$-compact in $PM_\Psi(G),$ to prove that $G_\mathbf{L}$ is compact, it is enough to prove that $0$ is not a $w^\ast$-cluster point of $\kappa(G_\mathbf{L}).$ Let us suppose to the contrary that $0$ is a $w^\ast$-cluster point of $\kappa(G_\mathbf{L}),$ i.e., there exists a net $(g_\alpha)$ in $G_\mathbf{L}$ such that $\lambda_\Psi(g_\alpha) \to 0.$ This implies that 
        $$u(g_\alpha) \to 0 \, \, \forall \, \, u \in \mathbf{L}.$$
        Observe that for $g \in G$ and $u \in \mathbf{L},$ since $g_\alpha \in G_\mathbf{L}$ for each $\alpha$ and $\mathbf{L}$ is right translation invariant, we have,
        $$u(g) = (L_{g_\alpha^{-1}}u)(g) = u(g_\alpha \, g) = (R_g u)(g_\alpha) \to 0.$$
        As $g \in G$ and $u \in \mathbf{L}$ is arbitrary, this implies that $\mathbf{L} = \{0\}$ which is a contradiction. Thus, the subgroup $G_\mathbf{L}$ is compact. Now, for this compact subgroup $G_\mathbf{L}$ of $G,$ consider 
        $$\mathcal{L}_\Phi(G_\mathbf{L}) = \{u \in A_\Phi(G): L_g u = u \, \, \forall \, \, g \in G_\mathbf{L}\}.$$ 
        It is clear that $\mathbf{L} \subseteq \mathcal{L}_\Phi(G_\mathbf{L}).$ To prove the equality, let $\overline{\mathbf{L}}^{w^\ast}$ be the $w^\ast$-closure of $\mathbf{L}$ in $L^\infty(G).$ Since $\mathbf{L}$ is a non-zero right translation invariant subalgebra of $A_\Phi(G) \subseteq L^\infty(G)$ which is closed under conjugation, $\overline{\mathbf{L}}^{w^\ast}$ is non-zero right translation invariant von Neumann subalgebra of $L^\infty(G).$ By \cite[Theorem 2]{tak}, there exists a unique closed subgroup $H_\mathbf{L}$ of $G$ such that 
        $$\overline{\mathbf{L}}^{w^\ast} = \{u \in L^\infty(G): L_h u = u \, \, \forall \, \, h \in H_\mathbf{L} \}$$
        and 
        $$H_\mathbf{L} = \{h \in G: L_h u = u \, \, \forall \, \, u \in \overline{\mathbf{L}}^{w^\ast}\}.$$
        Observe that $H_\mathbf{L} = G_\mathbf{L}$ and $\mathcal{L}_\Phi(G_\mathbf{L}) \subseteq \overline{\mathbf{L}}^{w^\ast}.$ Further, it is easy to verify that $\overline{\mathbf{L}}^{w^\ast} \cap A_\Phi(G) = \mathbf{L},$ which implies $\mathcal{L}_\Phi(G_\mathbf{L}) \subseteq \mathbf{L}$ and hence, $\mathcal{L}_\Phi(G_\mathbf{L}) = \mathbf{L}.$ Thus, it follows that $\mathcal{F}(G_\mathbf{L}) = \mathbf{L},$ proving that the map $\mathcal{F}$ is surjective.
            \end{proof}
            
        \begin{cor}
            Let $H$ be a closed subgroup of $G.$ Then $\mathcal{L}_\Phi(H) \neq \{0\}$ if and only if $H$ is compact.
        \end{cor}

        \begin{proof}
            By Lemma \ref{5.1}, if $H$ is compact subgroup of $G,$ then $\mathcal{L}_\Phi(H) \neq \{0\}.$

            Conversely, assume that $\mathcal{L}_\Phi(H) \neq \{0\}.$ By Theorem \ref{5.3}, there exists a unique compact subgroup, say $K,$ of $G$ such that $\mathcal{L}_\Phi(K) = \mathcal{L}_\Phi(H)$ and $$K = \{g \in G: L_g u = u \, \, \forall \, \, u \in \mathcal{L}_\Phi(H)\}.$$
            It is easy to verify that $H \subseteq K$ and hence, the result follows.
        \end{proof}

   \section*{Acknowledgement}
    The first author is grateful to the Indian Institute of Technology Delhi for the Institute Assistantship.

    \section*{Data Availability} 
    Data sharing does not apply to this article as no datasets were generated or analysed during the current study.

    \section*{Competing Interests}
    The authors declare that they have no competing interests.

	\bibliographystyle{acm}
	\bibliography{article4}

\end{document}